\pdfoutput=1
\documentclass[11pt]{amsart}

\usepackage{amsmath,amssymb,amsthm}
\usepackage[margin=1.2in]{geometry}
\usepackage{microtype}
\usepackage{booktabs}
\usepackage{needspace}
\usepackage{algorithm}
\usepackage[noend]{algpseudocode}
\usepackage{tikz}
\usepackage{xcolor}
\usepackage[colorlinks=true,linkcolor=blue!50!black,citecolor=blue!50!black,urlcolor=blue!50!black]{hyperref}
\usepackage{xurl}
\hypersetup{pdftitle={The square achievement game on n x n boards},
  pdfauthor={Sean Reid}}

\newtheorem{theorem}{Theorem}[section]
\newtheorem{lemma}[theorem]{Lemma}
\newtheorem{proposition}[theorem]{Proposition}
\theoremstyle{definition}
\newtheorem{definition}[theorem]{Definition}
\theoremstyle{remark}
\newtheorem{remark}[theorem]{Remark}

\newcommand{\Board}{B}
\newcommand{\Sq}{Q}

\newcommand{\pairsSix}{1{,}155}
\newcommand{\pairsSeven}{2{,}208}
\newcommand{\pairsEight}{3{,}843}
\newcommand{\pairsNine}{6{,}240}
\newcommand{\pairsTen}{9{,}603}
\newcommand{\pairsEleven}{14{,}160}
\newcommand{\pairsTwelve}{20{,}163}
\newcommand{\pairsThirteen}{27{,}888}
\newcommand{\pairsFourteen}{37{,}635}
\newcommand{\pairsFifteen}{49{,}728}
\newcommand{\pairsSixteen}{64{,}515}
\newcommand{\pairsSeventeen}{82{,}368}
\newcommand{\tempoZeroSeventeen}{81{,}186}
\newcommand{\tempoOneSeventeen}{1{,}103}
\newcommand{\pairsAnchor}{277{,}728}
\newcommand{\anchorAbstractions}{43{,}424}
\newcommand{\blocksThatThreaten}{19}
\newcommand{\distinctLines}{79}
\newcommand{\tempoZero}{985}
\newcommand{\tempoOne}{91}
\newcommand{\tempoTwo}{54}
\newcommand{\tempoThree}{25}
\newcommand{\longestLineMoves}{6}
\newcommand{\latestPly}{17}
\newcommand{\anchoredPositions}{176}
\newcommand{\fiveTablePositions}{1.79 \times 10^8}
\newcommand{\fiveTableBranches}{24}
\newcommand{\fiveNoForcing}{104}
\newcommand{\ruleLeanLines}{250}
\newcommand{\timeFive}{1.5}
\newcommand{\timeSix}{0.9}
\newcommand{\timeSeven}{1.3}
\newcommand{\timeEight}{1.9}
\newcommand{\timeNine}{2.9}
\newcommand{\timeTen}{4.4}
\newcommand{\timeEleven}{6.8}
\newcommand{\timeTwelve}{10}
\newcommand{\timeThirteen}{16}
\newcommand{\timeFourteen}{24}
\newcommand{\timeFifteen}{36}
\newcommand{\timeSixteen}{52}
\newcommand{\timeSeventeen}{72}
\newcommand{\timeAnchor}{499}
\newcommand{\timeAnchored}{0.34}
\newcommand{\buildMinutes}{17}

\title{The square achievement game on \texorpdfstring{$n \times n$}{n x n} boards}
\author{Sean Reid}
\address{Department of Mathematics and Statistics, San Jose State
  University, San Jose, CA 95192}
\email{sean.reid@sjsu.edu}

\begin{document}

\begin{abstract}
Two players alternately claim cells of an $n \times n$ grid; the first
to hold all four corners of an axis-parallel square wins, and a full
board without one is a draw. The first player has a winning strategy
if and only if $n \ge 5$, which answers a question of Erickson. For
$n \ge 6$ one strategy serves every board: a rule with no lookup table
whose moves after the second reply are all forcing, built around the
\emph{plus}, a stone with the four cells at a common distance from it
along the axes, a shape proposed for this game by Sirgedas in 2010. The
rule is checked on the boards of size $6$ to $17$, and a simulation
theorem carries one checked run on the $23 \times 23$ board to every
larger board. The proofs are formalized in Lean~4, and the theorem for
the rule states that the tree of its own moves is winning.
\end{abstract}

\maketitle

\section{Introduction}

Two players alternately claim empty cells of an $n \times n$ grid, the
first player moving first. A player wins upon holding all four corners
of a square with horizontal and vertical sides, of any size: the cells
$(r,c)$, $(r,c+s)$, $(r+s,c)$, $(r+s,c+s)$ for some $s \ge 1$. If the
board fills with no player holding such a square, the game is a draw.

Erickson posed the game and asked for the outcome under optimal play
\cite{erickson-pearls, opg, erickson-mo}. It is an achievement game in the sense of
Harary \cite{harary}, played in the Maker--Maker convention: both
players pursue the same goal, and a player who loses the race may still
secure a draw. The convention goes back to the positional games of
Hales and Jewett \cite{hales-jewett}; a recent unified treatment of
achievement conventions, with the strategy-stealing theory in the
generality used here, is \cite{galliot-senizergues}, and the general
theory of positional games is in Beck \cite{beck}.

Two boundary results were known. Jenrich \cite{jenrich} reported
computer searches showing that the second player can secure a draw for
$n = 3, 4$ and that the first player can force a win for $n = 5$.
Bacher and Eliahou \cite{bacher-eliahou} proved that every two-coloring
of the $15 \times 15$ grid contains four equal-colored cells at the
corners of an axis-parallel square, and that $15$ is least with this
property. Since no completed game on a board of size $15$ or more is
square-free, strategy stealing and determinacy give the first player a
win there. Between these results nothing was published: for
$6 \le n \le 14$ neither a winning strategy nor a drawing strategy was
known \cite{opg, adler}, and no coloring argument can reach this range,
since Bacher and Eliahou also exhibit square-free two-colorings of the
$14 \times 14$ board. One idea had been sketched. Answering Erickson's
question on MathOverflow in 2010, Sirgedas \cite{sirgedas} proposed
that the first player open in the center and win by making a
\emph{cross}, a stone with its four neighbors at a common distance
along the axes, each arm forcing a reply; he noted that the cross can
also be built on two stones in a row or column at distance two, worked
out three of the nine first replies on the $7 \times 7$ board up to
symmetry, and guessed that $n = 5$ was a draw. Concurrently with the
present work, Adler \cite{adler, adler-code} announced computer
searches settling the nine sizes $6 \le n \le 14$ as first-player
wins, with full strategy trees certified for $n = 6$ and $7$ and search
verdicts for the rest, cross-checked by independent runs at the smaller
sizes; Adler notes that the results are unrefereed and names a
conceptual proof with an explicit strategy as the open direction.

\begin{theorem}\label{thm:main}
In the square achievement game on the $n \times n$ board, the first
player has a winning strategy if and only if $n \ge 5$. For $n \le 4$
the game is drawn under optimal play: each player can secure at least a
draw.
\end{theorem}

Relative to Adler's announcement, what is new is a single strategy
that wins every board from $6 \times 6$ up, stated in a paragraph, and
a proof checked by a proof assistant. The strategy is a rule built
around one shape, the plus of Definition~\ref{def:plus}, which is
Sirgedas's cross: the arms of a plus force the second player's replies
one at a time, and the last arm threatens two squares at once. The
sketch of 2010 covered one board and a few replies and did not say
which plus to aim for, how to prepare one when the replies obstruct
it, or why one board's analysis should carry to another. The rule
settles the first two by a search with a fixed preference order, and
the simulation theorem of Section~\ref{sec:uniform} settles the
third.

\begin{theorem}\label{thm:rule}
For every $n \ge 6$, the first player following the plus rule of
Section~\ref{sec:rule} wins the square achievement game on the
$n \times n$ board.
\end{theorem}

Two facts make Theorem~\ref{thm:rule} a finite computation. After the
second reply the rule makes only forcing moves, so every later reply
is a forced block and the game tree branches freely at two levels only.
Checking the tree on one board means enumerating the second player's
first two replies and following a chain of forced moves from each pair;
this check is run on every board from $6$ to $17$
(Section~\ref{sec:verify}). The rule reads only the top-left
$17 \times 17$ sub-board and one cell computed from the second player's
first two replies. A stone outside the sub-board can affect the game
only through that cell, and the rule never forces a block there, so a
game on any larger board is mirrored move for move by a game on the
$23 \times 23$ board. A simulation theorem, stated for any strategy, turns one checked
run on that board into a proof for every board of size at least $17$
(Section~\ref{sec:uniform}). Theorem~\ref{thm:main} then follows with an
exhaustive refutation for $n \le 4$ and a strategy table for $n = 5$,
whose play from the opening $(2,2)$ is not purely forcing
(Section~\ref{sec:small}).

Everything is formalized in Lean~4 \cite{lean4} on mathlib
\cite{mathlib} (Section~\ref{sec:lean}). The formal theorem for the
rule states that the tree of the rule's own moves is winning, not only
that a winning tree exists; the runs of the checkers are discharged by
evaluation, and the two properties of the rule that the simulation
theorem needs are proved as lemmas about its definition. The rule was
extracted from an earlier lookup-table strategy, almost all of whose
forcing sequences were arms of a side-$1$ plus at $(2,2)$ or $(3,2)$,
the rest preparing the side-$2$ plus at $(3,2)$. The $2011$ versions of
\cite{jenrich} contained an error in the search program, corrected in
$2012$; the formalization makes every search here discharge a
statement about the game rather than about a program.

\section{Definitions}\label{sec:defs}

Cells are pairs $(r,c)$ of nonnegative integers, row then column,
with row $0$ at the top and column $0$ at the left in the figures,
where the first player's stones are black and the second player's
white, and $\Board_n = \{0,\dots,n-1\}^2$. For $s \ge 1$ the square
with top-left corner $(r,c)$ and side $s$ is
$\Sq(r,c,s) = \{(r,c), (r,c+s), (r+s,c), (r+s,c+s)\}$; a set of cells
contains a square if it contains some $\Sq(r,c,s)$. The case $s = 1$ is
allowed, so the smallest winning shape is a $2 \times 2$ array, the
convention of \cite{bacher-eliahou, jenrich}; tilted squares do not
count. A \emph{position} is a pair of disjoint subsets of $\Board_n$,
the stones of the two players. The game is finite with perfect
information, so exactly one of three holds: the first player can force
a win, the second can, or both can force at least a draw.

A player \emph{threatens} at a cell $z$ if $z$ is empty and $z$
together with three of that player's stones is a square; the opponent
\emph{blocks} by playing $z$. A move that makes exactly one threat is
\emph{forcing}; a move that makes two threats wins if the opponent has
no threat of their own, since one block leaves the other. A \emph{ply} is one move of either player,
so the first player moves at the odd plies. The rule of
Section~\ref{sec:rule} plays inside the \emph{region} $\Board_6$ and
reads the second player's stones inside the \emph{base board}
$\Board_{17}$.

\subsection{Pluses}

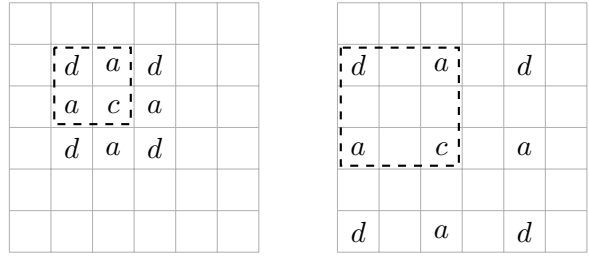
\begin{figure}[ht]
\centering
% Generated by python3 -m squares.paper.figures; do not edit.
\begin{tikzpicture}[scale=0.55]
  \draw[gray!60] (0,0) grid (6,-6);
  \node at (2.5,-2.5) {$c$};
  \foreach \r/\c in {1/2, 3/2, 2/1, 2/3} \node at (\c+0.5,-\r-0.5) {$a$};
  \foreach \r/\c in {1/1, 1/3, 3/1, 3/3} \node at (\c+0.5,-\r-0.5) {$d$};
  \draw[thick,dashed] (1.08,-1.08) rectangle (2.92,-2.92);
\end{tikzpicture}
\hspace{2em}
\begin{tikzpicture}[scale=0.55]
  \draw[gray!60] (0,0) grid (6,-6);
  \node at (2.5,-3.5) {$c$};
  \foreach \r/\c in {1/2, 5/2, 3/0, 3/4} \node at (\c+0.5,-\r-0.5) {$a$};
  \foreach \r/\c in {1/0, 1/4, 5/0, 5/4} \node at (\c+0.5,-\r-0.5) {$d$};
  \draw[thick,dashed] (0.08,-1.08) rectangle (2.92,-3.92);
\end{tikzpicture}
\caption{Left: the plus of side $1$ at $c = (2,2)$, with its arms $a$
and diagonals $d$. Right: the plus of side $2$ at $c = (3,2)$. In each,
the dashed outline is one of the four squares of the plus, whose
corners are the center, two perpendicular arms, and the diagonal
between them.}
\label{fig:anatomy}
\end{figure}

\begin{definition}\label{def:plus}
The \emph{plus} of side $s$ at the \emph{center} $c = (r,k)$ has four
\emph{arms}, the cells $(r-s,k)$, $(r+s,k)$, $(r,k-s)$, $(r,k+s)$, and
four \emph{diagonals}, the cells $(r \pm s, k \pm s)$
(Figure~\ref{fig:anatomy}). Two arms are \emph{perpendicular} if one
is above or below the center and the other to its left or right; the
diagonal between them is the cell in the row of one and the column of
the other. A plus lies in the region if its arms and diagonals do, that
is, $s \le r, k$ and $r+s, k+s \le 5$.
\end{definition}

The region holds $16$ pluses of side $1$, with centers in
$\{1,\dots,4\}^2$, and four of side $2$, with centers in $\{2,3\}^2$;
side $3$ does not fit.

\begin{lemma}\label{lem:arms}
Suppose the first player holds the center $c$ and one arm $a$ of a plus
of side $s$ whose four diagonals are empty. Playing an empty arm $b$
perpendicular
to $a$ threatens at the
diagonal between $a$ and $b$. If the first player also holds the arm
$a'$ opposite to $a$, playing $b$ threatens at two cells, the diagonals
between $a$ and $b$ and between $a'$ and $b$.
\end{lemma}

\begin{proof}
Say $a = (r-s,k)$ and $b = (r,k-s)$. Then $c$, $a$, $b$ and
$d = (r-s,k-s)$ are the corners of $\Sq(r-s,k-s,s)$, so with $c, a, b$
held and $d$ empty the first player threatens at $d$. With
$a' = (r+s,k)$ also held, $c$, $a'$, $b$ and $(r+s,k-s)$ are the
corners of $\Sq(r,k-s,s)$, and $(r+s,k-s) \ne d$. The other placements
are symmetric.
\end{proof}

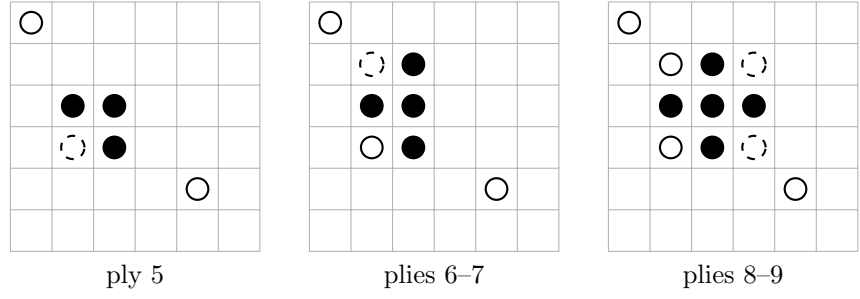
\begin{figure}[ht]
\centering
% Generated by python3 -m squares.paper.figures; do not edit.
\begin{tikzpicture}[scale=0.55]
  \draw[gray!60] (0,0) grid (6,-6);
  \foreach \r/\c in {2/2, 3/2, 2/1}
    \filldraw[fill=black,thick] (\c+0.5,-\r-0.5) circle (0.26);
  \foreach \r/\c in {0/0, 4/4}
    \filldraw[fill=white,thick] (\c+0.5,-\r-0.5) circle (0.26);
  \draw[thick,dashed] (1.5,-3.5) circle (0.28);
  \node at (3,-6.6) {\small ply 5};
\end{tikzpicture}
\hspace{1em}
\begin{tikzpicture}[scale=0.55]
  \draw[gray!60] (0,0) grid (6,-6);
  \foreach \r/\c in {2/2, 3/2, 2/1, 1/2}
    \filldraw[fill=black,thick] (\c+0.5,-\r-0.5) circle (0.26);
  \foreach \r/\c in {0/0, 4/4, 3/1}
    \filldraw[fill=white,thick] (\c+0.5,-\r-0.5) circle (0.26);
  \draw[thick,dashed] (1.5,-1.5) circle (0.28);
  \node at (3,-6.6) {\small plies 6--7};
\end{tikzpicture}
\hspace{1em}
\begin{tikzpicture}[scale=0.55]
  \draw[gray!60] (0,0) grid (6,-6);
  \foreach \r/\c in {2/2, 3/2, 2/1, 1/2, 2/3}
    \filldraw[fill=black,thick] (\c+0.5,-\r-0.5) circle (0.26);
  \foreach \r/\c in {0/0, 4/4, 3/1, 1/1}
    \filldraw[fill=white,thick] (\c+0.5,-\r-0.5) circle (0.26);
  \draw[thick,dashed] (3.5,-1.5) circle (0.28);
  \draw[thick,dashed] (3.5,-3.5) circle (0.28);
  \node at (3,-6.6) {\small plies 8--9};
\end{tikzpicture}
\caption{The basic line (Definition~\ref{def:line}) against the replies
$(0,0)$ and $(4,4)$. The first player holds $(2,2)$ and $(3,2)$; each
arm threatens at the dashed cell, the second player blocks, and the
last arm threatens twice.}
\label{fig:basic}
\end{figure}

So from a center and one arm the first player reaches a double threat
in three moves against a second player who has no threat and whose
blocks make none.
The two arms perpendicular to the held arm each force a block on a
diagonal, and the opposite arm threatens at two diagonals at once
(Figure~\ref{fig:basic}). The opening $(2,2)$ and the second move
$(3,2)$ are the center and an arm of the side-$1$ plus at $(2,2)$, and
an arm and the center of the side-$1$ plus at $(3,2)$.

\subsection{Far stones and the hidden corner}

The rule's stones never leave $\Board_6$, and the second player's
stones in a position the rule reaches are two free replies and forced
blocks, the blocks in $\Board_6$ as well.

\begin{lemma}\label{lem:span}
Let $x \ne y$ be distinct corners of a square $\Sq(r,c,s)$ with
$x, y \in \Board_m$. Then $\Sq(r,c,s) \subseteq \Board_{2m}$, and every
corner of $\Sq(r,c,s)$ has its row or its column indexed below $m$.
\end{lemma}

\begin{proof}
Two distinct corners either share a coordinate and differ by $s$ in the
other, or differ by $s$ in both; either way $s \le m - 1$. Each
coordinate of $(r,c)$ is bounded by the corresponding coordinate of a
corner, so $r, c \le m - 1$, and the containment follows, indeed with
$r + s, c + s \le 2m - 2$. If $x$ and $y$ share a row,
the two remaining corners inherit their columns, both below $m$; a
shared column is symmetric; and if $x$ and $y$ are diagonally opposite,
all four corners lie in $\Board_m$.
\end{proof}

\begin{figure}[ht]
\centering
% Generated by python3 -m squares.paper.figures; do not edit.
\begin{tikzpicture}[scale=0.3]
  \fill[gray!25] (0,0) rectangle (6,-6);
  \draw[gray!60] (0,0) grid (18,-18);
  \draw[semithick] (0,0) rectangle (17,-17);
  \draw[very thick] (0,0) rectangle (18,-18);
  \draw[thick,dashed] (3.5,-0.5) rectangle (17.5,-14.5);
  \foreach \r/\c in {0/17, 14/3}
    \filldraw[fill=white,thick] (\c+0.5,-\r-0.5) circle (0.42);
  \draw[thick,dashed] (3.5,-0.5) circle (0.42);
  \draw[thick,dashed] (17.5,-14.5) circle (0.42);
  \node[anchor=south,inner sep=1pt] at (3.5,0.04) {$b$};
  \node[anchor=south,inner sep=1pt] at (17.5,0.04) {$o_1$};
  \node[anchor=east,inner sep=1pt] at (2.96,-14.5) {$o_2$};
  \node[anchor=west,inner sep=1pt] at (18.04,-14.5) {$z$};
\end{tikzpicture}
\caption{A hidden corner, on $\Board_{18}$. The inner outline is
$\Board_{17}$ and the shaded square the region $\Board_6$. The second
player's free replies
$o_1 = (0,17)$ and $o_2 = (14,3)$ are two corners of the dashed
square; its third corner $b = (0,3)$ is their hidden corner, and $o_1$
and the fourth corner $z = (14,17)$ lie beyond $\Board_{17}$. A second
player forced to block at $b$ would threaten at $z$, a cell the rule
does not read.}
\label{fig:latent}
\end{figure}
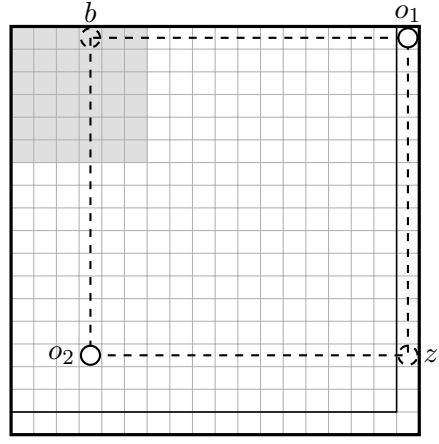

With $m = 6$, two stones of the region lie in a common square only if
that square lies in $\Board_{12}$. Hence a second-player stone beyond
$\Board_{17}$ shares no square with two blocks. A second-player threat
through such a stone and a block must therefore have the other free
reply as its third corner (Figure~\ref{fig:latent}). The
block in that situation is the one cell through which stones beyond
$\Board_{17}$ can affect the game.

\begin{definition}\label{def:latent}
Let $x$ and $y$ be two cells each having a coordinate at least $12$.
Their \emph{hidden corner} is a cell $b \in \Board_6$ such that
$b$, $x$, $y$ are three corners of one square. The \emph{abstraction}
of a position consists of the first player's stones, the second
player's stones inside $\Board_{17}$, and, once the second player has
two stones, the hidden corner of the first two, if they have one.
\end{definition}

\begin{lemma}\label{lem:dichotomy}
Two cells $x \ne y$ have at most one hidden corner. If $\Sq(r,c,s)$ has
a corner $b \in \Board_6$, has $x$ and $y$ among its other corners, and
is not contained in $\Board_{17}$, then $b = (r,c)$, $s \ge 12$, each
of $x$, $y$ has a coordinate at least $12$, and $b$ is the hidden
corner of $\{x, y\}$.
\end{lemma}

\begin{proof}
Uniqueness. If $x$ and $y$ share a column, the square has side
$s = |x_1 - y_1|$, and $b$ lies in the row of a member with row below
$6$, at distance $s$ to its left or right. Only one member can have row
below $6$: if both did, then $s < 6$ and the shared column, being the
coordinate at least $12$ of either, would put both candidates for $b$
at column $7$ or more. For that member the cell to the right is
excluded too: it lies in $\Board_6$ only if the member's column is
below $6$, and that column is the member's coordinate at least $12$,
its row being below $6$. A shared row is symmetric. If $x$ and $y$ are diagonally
opposite, the two remaining corners are $(x_1, y_2)$ and $(y_1, x_2)$,
of which at most one lies in $\Board_6$, since otherwise all four
coordinates of $x$ and $y$ would be below $6$.

For the second claim, a second corner in $\Board_6$ would put the
square inside $\Board_{12}$ by Lemma~\ref{lem:span}, against the
hypothesis; so $b$ is alone in the region. Were $b$ not $(r,c)$, the
corner $(r,c)$ would have both coordinates at most those of $b$ and
would lie in $\Board_6$; so $b = (r,c)$. The corner opposite $b$ has
both coordinates at most $5 + s$, and the square leaves $\Board_{17}$
only if a coordinate reaches $17$, so $s \ge 12$. Each of $x$, $y$ has
a coordinate equal to $r + s$ or $c + s$, hence at least $12$, and $b$
is their hidden corner by definition.
\end{proof}

The rule therefore never makes a threat whose only block is the hidden
corner (Definitions~\ref{def:tempo} and~\ref{def:line}). It applies the
restriction whenever the pair has a hidden corner, even when the square
through it lies inside $\Board_{17}$ and the restriction is not
needed.

\begin{remark}\label{rem:base}
The threshold $12$ of Definition~\ref{def:latent} and the base $17$ are
tied by $17 = 12 + 5$: a square through a region cell leaves
$\Board_{17}$ only if its side is at least $12$. With a threshold $t$
and the base $t + 5$ in their place, Lemma~\ref{lem:dichotomy} holds
exactly when $t \ge 11$. At $t = 10$ uniqueness fails: $(0,10)$ and
$(5,10)$ would have the two hidden corners $(0,5)$ and $(5,5)$. So base
$16$ with threshold $11$ would serve the lemma as well. The development
uses $t = 2 \cdot 6 = 12$, the bound of Lemma~\ref{lem:span}, and the
rule has been checked only with base $17$.
\end{remark}

\subsection{Tempo moves and lines}

When no plus can be completed directly, because the second player holds
an arm or a diagonal or because the arms would not force, the rule
first plays moves that each force a block and so add a stone. Threats
in this subsection are those of the abstraction: squares of the first
player inside $\Board_6$ and squares of the second player inside
$\Board_{17}$.

\begin{definition}\label{def:tempo}
A \emph{tempo move} in a position is a cell of $\Board_6$ such that,
after it is played, the second player has no threat and the first
player has exactly one, at a cell other than the hidden corner, if the
pair has one.
\end{definition}

In a position where the second player has a threat, a tempo move must
remove it, and the only way is to occupy its cell; so the second player
has exactly one threat there and the tempo move is its block.

\begin{definition}\label{def:line}
A \emph{line} from a position consists of some tempo moves, each
followed by its block, and then, in some order, the arms the first
player does not yet hold of one plus in the region whose center the
first player holds and none of whose arms the second player holds. The
arms are played until one makes at least two threats. After each move
of the line the second player has no threat, and each arm before the
last played makes exactly one threat, at a cell other than the hidden
corner. The \emph{key} of a
line is
\[
  (\text{number of tempo moves},\ \text{number of arms not held},\
  \text{side of the plus},\ \text{center of the plus}).
\]
The rule prefers the smallest key in lexicographic order, cells ordered
by row and then by column, then the line its search generates first.
\end{definition}

The arms after the one that makes two threats are part of the line and
of its key but are never played. A block may give the second player a
threat; the next move of the line is then a tempo move that takes it.
The rule searches for lines with at most three tempo moves. Whether a
line exists, and which is preferred, depends only on the abstraction of
the position.

\section{The plus rule}\label{sec:rule}

\begin{figure}[ht]
\centering
% Generated by python3 -m squares.paper.figures; do not edit.
\begin{tikzpicture}[scale=0.55]
  \draw[gray!60] (0,0) grid (6,-6);
  \foreach \r/\c in {2/2}
    \filldraw[fill=black,thick] (\c+0.5,-\r-0.5) circle (0.26);
  \foreach \r/\c/\t in {1/3/1, 2/1/2, 2/3/1, 3/1/3, 3/2/3, 3/3/1, 4/3/2}
    \node at (\c+0.5,-\r-0.5) {$\t$};
\end{tikzpicture}
\caption{The second move. The black stone is the opening. Against a
first reply on a cell marked $1$ the second move is $(2,1)$; on $2$ it
is $(2,3)$; on $3$ it is $(1,2)$. Against every other first reply, on
this board or beyond it, it is $(3,2)$.}
\label{fig:second}
\end{figure}
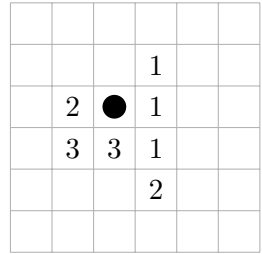

Algorithm~\ref{alg:move} is the move function and
Algorithm~\ref{alg:line} the line search it calls once, after the
second reply. Both are stated as in the formal development. The second
move $x_2(o_1)$ against a first reply $o_1$ is $(3,2)$ unless $o_1$ is
one of seven cells, for which Figure~\ref{fig:second} gives the move.

\begin{algorithm}[H]
\caption{The move function.}\label{alg:move}
\begin{algorithmic}[1]
\If{some empty cell completes a square of the first player}
  \State play the least such cell
\ElsIf{some empty cell completes a square of the second player}
  \State play the least such cell
\ElsIf{the first player has no stone}
  \State play $(2,2)$
\ElsIf{the first player does not hold $x_2(o_1)$}
  \State play $x_2(o_1)$
\Else
  \State $L \gets$ the line of Algorithm~\ref{alg:line} from the
    abstraction after the second reply
  \State play the first cell of $L$ the first player does not hold
\EndIf
\end{algorithmic}
\end{algorithm}

Steps 1 to 4 scan the squares of the abstraction: the first player's,
which lie in $\Board_6$, and the second player's with all corners
inside $\Board_{17}$. A second-player square with a corner beyond
$\Board_{17}$ would pass through the hidden corner, where the rule never
forces a block. Independently, the checker of Section~\ref{sec:verify}
confirms at every node that the second player has no threat at all. If
Algorithm~\ref{alg:line} returns no line the rule is undefined; the
checker certifies that this never happens at a position the rule
reaches.

\begin{algorithm}[H]
\caption{The line search, from an abstract position $A$.}\label{alg:line}
\begin{algorithmic}[1]
\State $F \gets [(A, \varepsilon)]$
  \Comment{positions reached, with the tempo moves that reached them}
\For{$d = 0, 1, 2, 3$}
  \State $R \gets \varnothing$
  \ForAll{$(X, p) \in F$ in order, such that the second player has no threat in $X$}
    \ForAll{pluses $P$ in the region whose center the first player holds in $X$}
      \State $T \gets$ the arms of $P$ the first player does not hold in $X$
      \If{$T \ne \varnothing$ and the second player holds no cell of $T$}
        \State $\pi \gets$ the first ordering of $T$ such that $p\,\pi$ is a line from $A$
        \If{$\pi$ exists} record $p\,\pi$ in $R$ with the key $(d, |T|, \text{side of } P, \text{center of } P)$ \EndIf
      \EndIf
    \EndFor
  \EndFor
  \If{$R \ne \varnothing$} \Return the first line recorded in $R$ among those with the smallest key \EndIf
  \State $F \gets$ the successors of the members of $F$, in order
\EndFor
\State \Return no line
\end{algorithmic}
\end{algorithm}

The successors of $(X, p)$ in step 11 are the pairs $(X', p\,m)$ for
the tempo moves $m$ in $X$, taken in cell order, where $X'$ is $X$
after $m$ and its block. They are taken from every member of $F$,
including positions in which the second player has a threat; such a
position offers no plus (step 4), but its tempo moves, which must
remove the threat, are explored. Step 8 tries the orderings of $T$ in
lexicographic order, the arms listed above, below, left of, and right
of the center, so that two implementations choose the same line.

The exceptions are exactly the first replies against which $(3,2)$ is
occupied or fails. They were found by search: for each first reply on $\Board_6$ the
second move is the first cell in the order $(3,2)$, $(2,3)$, $(2,1)$,
$(1,2)$ against which the rest of the rule beats every second reply.

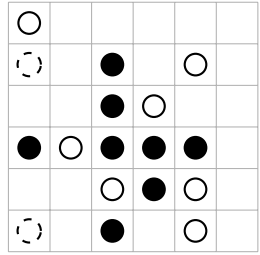
\begin{figure}[ht]
\centering
% Generated by python3 -m squares.paper.figures; do not edit.
\begin{tikzpicture}[scale=0.55]
  \draw[gray!60] (0,0) grid (6,-6);
  \foreach \r/\c in {2/2, 3/2, 3/3, 4/3, 3/4, 1/2, 5/2, 3/0}
    \filldraw[fill=black,thick] (\c+0.5,-\r-0.5) circle (0.26);
  \foreach \r/\c in {0/0, 3/1, 2/3, 4/2, 4/4, 1/4, 5/4}
    \filldraw[fill=white,thick] (\c+0.5,-\r-0.5) circle (0.26);
  \draw[thick,dashed] (0.5,-1.5) circle (0.28);
  \draw[thick,dashed] (0.5,-5.5) circle (0.28);
\end{tikzpicture}
\caption{The second game of Table~\ref{tab:games} after ply $15$. The
side-$2$ plus at $(3,2)$ stands with the second player's blocks on its
diagonals $(1,4)$ and $(5,4)$, and the last arm $(3,0)$ threatens at
the other two diagonals, dashed.}
\label{fig:longest}
\end{figure}

Table~\ref{tab:games} lists two complete games in which the second
player always blocks. In the first the replies $(0,0)$ and $(4,4)$ are
away from the opening pair and the basic line wins at once, key
$(0,3,1,(2,2))$.
In the second the reply $(3,1)$ is a diagonal of the side-$1$ plus at
$(2,2)$ and an arm of the side-$1$ plus at $(3,2)$,
no line with fewer than two tempo moves exists, and the line has key
$(2,4,2,(3,2))$: two tempo moves and then the four arms of the side-$2$
plus at $(3,2)$. Figure~\ref{fig:longest} shows the position before its
last move. After the forced block at ply $10$ of the first game and at
ply $16$ of the second, the second player does hold a threat, at the
cell the first player completes next; Definition~\ref{def:line}
constrains the second player's threats only after the first player's
moves.

\begin{table}[ht]
\centering
\small
\begin{tabular}{rllp{0.5\textwidth}}
\toprule
ply & first player & second player & what the move does \\
\midrule
1, 2 & $(2,2)$ & $(0,0)$ & opening; free reply \\
3, 4 & $(3,2)$ & $(4,4)$ & second move; free reply \\
5, 6 & $(2,1)$ & $(3,1)$ & arm; threat at $(3,1)$ through $\Sq(2,1,1)$, blocked \\
7, 8 & $(1,2)$ & $(1,1)$ & arm; threat at $(1,1)$ through $\Sq(1,1,1)$, blocked \\
9, 10 & $(2,3)$ & $(1,3)$ & arm; threats at $(1,3)$ and $(3,3)$; one is blocked \\
11 & $(3,3)$ & & completes $\Sq(2,2,1)$; wins \\
\midrule
1, 2 & $(2,2)$ & $(0,0)$ & opening; free reply \\
3, 4 & $(3,2)$ & $(3,1)$ & second move; free reply on a diagonal \\
5, 6 & $(3,3)$ & $(2,3)$ & tempo move; threat at $(2,3)$ through $\Sq(2,2,1)$, blocked \\
7, 8 & $(4,3)$ & $(4,2)$ & tempo move; threat at $(4,2)$ through $\Sq(3,2,1)$, blocked \\
9, 10 & $(3,4)$ & $(4,4)$ & arm; threat at $(4,4)$ through $\Sq(3,3,1)$, blocked \\
11, 12 & $(1,2)$ & $(1,4)$ & arm; threat at $(1,4)$ through $\Sq(1,2,2)$, blocked \\
13, 14 & $(5,2)$ & $(5,4)$ & arm; threat at $(5,4)$ through $\Sq(3,2,2)$, blocked \\
15, 16 & $(3,0)$ & $(1,0)$ & arm; threats at $(1,0)$ and $(5,0)$; one is blocked \\
17 & $(5,0)$ & & completes $\Sq(3,0,2)$; wins \\
\bottomrule
\end{tabular}
\medskip
\caption{The basic line and the longest line as complete games.}
\label{tab:games}
\end{table}

The line is chosen once, after the second reply, and followed. Since
the blocks are forced, the line determines every position of the forced
tree (Section~\ref{sec:verify}) below the second reply. Until the
completing move, the first player's stones are the opening pair
together with a prefix of the line. Committing to the line makes the
forced tree a function of the reply pair, which
Section~\ref{sec:uniform} uses.

The Python implementation of the rule was run over every reply pair on
each board from $6$ to $17$, on the anchor board $\Board_{23}$, and
from every anchored position, both of Section~\ref{sec:uniform},
recording the line, its tempo moves, and the ply of the win when the
second player blocks. On every board from $6$ to $17$, on the anchor
board, and from every anchored position, the rule plays one of the same
$\distinctLines$ lines; the $\pairsAnchor$ pairs on the anchor board
have $\anchorAbstractions$ distinct abstractions. Every
abstraction that occurs on a board beyond $17$ occurs on the anchor
board or in an anchored position (Section~\ref{sec:uniform}), so the
same holds on every board. These counts come from the Python
implementation and the lines the checker prints; they are not part of
the formal theorem. Table~\ref{tab:lines} splits the pairs by
the number of tempo moves. Every pair won by a side-$1$ plus needs no
tempo move and every pair won by a side-$2$ plus needs at least one, so
the first column also counts the side-$1$ pairs. In
$\blocksThatThreaten$ of the $\pairsSix$ pairs on $\Board_6$, all with
three tempo moves, a block gives the second player a threat that the
next tempo move takes. The longest line, the second game of
Table~\ref{tab:games}, has $\longestLineMoves$ first-player moves after
the opening pair, and every game in which the second player blocks
ends by ply $\latestPly$; a reply that blocks nothing leaves a threat
that the rule completes at its next move. So the rule wins every game by the
$\latestPly$th ply. This bound is computed from the $\distinctLines$
lines and is not part of the formal theorem.

\begin{table}[ht]
\centering
\begin{tabular}{lrrrr}
\toprule
tempo moves in the line & 0 & 1 & 2 & 3 \\
\midrule
reply pairs on $\Board_6$ & \tempoZero & \tempoOne & \tempoTwo & \tempoThree \\
reply pairs on $\Board_{17}$ & \tempoZeroSeventeen & \tempoOneSeventeen & \tempoTwo & \tempoThree \\
\bottomrule
\end{tabular}
\medskip
\caption{The rule's lines by the number of tempo moves. The
$\tempoTwo + \tempoThree$ pairs needing two or three tempo moves have
both replies in $\Board_6$ and are the same on every board; the other
columns grow with the board.}
\label{tab:lines}
\end{table}

\section{Small boards and the \texorpdfstring{$5 \times 5$}{5 x 5} board}\label{sec:small}

\begin{lemma}\label{lem:stealing}
The second player has no winning strategy.
\end{lemma}

\begin{proof}
Suppose the second player had a winning strategy $\tau$. Let the first
player open with an arbitrary move and thereafter follow $\tau$ as if
the extra stone were absent, making a further arbitrary move whenever
$\tau$ selects an occupied cell. The winning condition is monotone in a
player's stone set, so surplus stones never hurt (a surplus stone may
itself complete a square, which ends the game in the first player's
favor), and the first player wins; but $\tau$ was assumed to win
against all first-player play, a contradiction.
\end{proof}

\begin{proposition}\label{prop:small}
For $n \le 4$ the first player has no winning strategy. With
Lemma~\ref{lem:stealing}, these games are drawn under optimal play:
each player can secure at least a draw.
\end{proposition}

\begin{proof}
For $n \le 2$ the board has at most four cells, so neither player ever
holds four. For $n = 3, 4$ the claim is a statement about a finite game
tree, verified by an exhaustive search that is part of the formal
development: for every first-player move the search exhibits a reply
after which the first player still cannot force a win. The result
agrees with \cite{jenrich}.
The step from ``neither player can force a win'' to ``each player can
secure a draw'' is determinacy of the finite game, proved by induction
on the number of empty cells.
\end{proof}

\Needspace*{6\baselineskip}
\begin{proposition}\label{prop:five}
The first player has a winning strategy on $\Board_5$.
\end{proposition}

\begin{proof}
The strategy is a table, verified by a checker that walks its complete
game tree on $\Board_5$, every reply at every node, and finds every
branch ending in a first-player square. The checker's soundness theorem
and the run are part of the formal development
(Section~\ref{sec:lean}).
\end{proof}

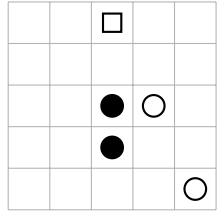
\begin{figure}[ht]
\centering
% Generated by python3 -m squares.paper.figures; do not edit.
\begin{tikzpicture}[scale=0.55]
  \draw[gray!60] (0,0) grid (5,-5);
  \foreach \r/\c in {2/2, 3/2}
    \filldraw[fill=black,thick] (\c+0.5,-\r-0.5) circle (0.26);
  \foreach \r/\c in {2/3, 4/4}
    \filldraw[fill=white,thick] (\c+0.5,-\r-0.5) circle (0.26);
  \draw[thick] (2.28,-0.72) rectangle (2.72,-0.28);
\end{tikzpicture}
\caption{A quiet move of the $5 \times 5$ table. The first player
opened at $(2,2)$ and answered the reply $(4,4)$ with the
table's second move $(3,2)$. After the second reply $(2,3)$, one of the
$\fiveNoForcing$ pairs with no forcing continuation, the table plays
$(0,2)$, marked by the square, which makes no threat, and then consults
a subtable indexed by the reply to it.}
\label{fig:five}
\end{figure}

The table was extracted from a full solution of the $5 \times 5$ game
($\fiveTablePositions$ positions solved, with draws scored against the
first player). Like the rule, it completes a square when it can and
otherwise blocks one before consulting the table. It opens at $(2,2)$,
gives a second move for each of the $\fiveTableBranches$ first replies,
and after the second reply plays a sequence of forcing moves ending in
a double threat where one exists. For $\fiveNoForcing$ of the reply
pairs, with the table's second move, the solver finds no winning
continuation of at most eight moves in which every first-player move
makes a threat, so the table plays a \emph{quiet} move there, one
making no threat, and continues with a subtable indexed by the reply
(Figure~\ref{fig:five}). That no such continuation exists is a
computation of the solver and is not formalized; the win of the table
is.

\section{Verification}\label{sec:verify}

Fix a strategy $\sigma$ and a board size $n$. The \emph{full tree} of
$\sigma$ on $\Board_n$ is the game tree in which the first player's
moves are given by $\sigma$ and the second player's moves range over
all empty cells; if every branch ends in a first-player square, the
first player wins. The
$5 \times 5$ table is checked on its full tree.

\Needspace*{10\baselineskip}
The rule is checked on a smaller tree. The \emph{forced tree} of
$\sigma$ on $\Board_n$ has the first player playing $\sigma$'s move,
the second player's replies ranging over all empty cells at the first
two levels, and below every later node only over the cells that block
a threat of the move just played. The checker walks this tree under
guards:
\begin{enumerate}
\item[(G1)] after the second reply, every move of $\sigma$ wins
  outright or makes a threat;
\item[(G2)] every nonwinning move and every threat cell lies in
  $\Board_6$;
\item[(G3)] after every move of $\sigma$ the second player has no
  threat;
\item[(G4)] optionally, for a given cell $g$, no move of $\sigma$ has
  $g$ as its only threat cell.
\end{enumerate}
A reply that blocks nothing leaves a threat standing, which the
strategy completes at its next move, so every omitted reply loses in
two plies.

\begin{proposition}\label{prop:forced}
Let $\sigma$ be a strategy that, whenever the first player's stones lie
in $\Board_6$ and some cell completes one of the first player's
squares, plays such a cell. If the forced tree of $\sigma$ on
$\Board_n$ passes the guards and every branch ends in a first-player
square, then the tree in which the first player plays $\sigma$'s move
at every node is winning; in particular the first player wins on
$\Board_n$.
\end{proposition}

\begin{proof}
This is the soundness theorem of the forced checker, an induction over
the checker's recursion (Section~\ref{sec:lean}). At a node after the
second reply, the move $m$ of $\sigma$ either completes a square or, by
(G1), makes a threat. By (G3) no reply completes a second-player
square. A reply that blocks a threat of $m$ is a child in the forced
tree and is winning by induction. A reply that blocks nothing leaves a
threat of $m$ standing, the first player's stones lie in $\Board_6$ by
(G2), and the hypothesis makes $\sigma$ complete the threat at its next
move, so that the winning tree is $\sigma$'s own.
\end{proof}

The rule satisfies the hypothesis: its first clause plays a completion
found by its threat scan, and that scan agrees with the checker's
whenever the first player's stones lie in $\Board_6$, which is proved
as a lemma (Section~\ref{sec:lean}).

The rule passes the forced tree on every board from $6$ to $17$, one
run of the checker per board with $(n^2-1)(n^2-3)$ reply pairs at the
root (Table~\ref{tab:runs} in Section~\ref{sec:lean}). A Python
implementation of the rule and a full-tree verifier were used while the
rule was being found. The checker executable prints the line the Lean
rule commits to for every reply pair of $\Board_6$ and for every
anchored position of Section~\ref{sec:uniform}, and a test of the
Python implementation reproduces every one; the full-tree verifier
confirms the win on $\Board_6$. Only the Lean runs enter the proof.

\section{Boards of size at least 17}\label{sec:uniform}

\begin{theorem}\label{thm:uniform}
For every $n \ge 17$ the first player following the rule wins on
$\Board_n$.
\end{theorem}

By Lemma~\ref{lem:dichotomy}, a stone beyond $\Board_{17}$ can affect
the rule's play only through a hidden corner, and the rule never
forces a block there. The proof replays a checked run on $\Board_{23}$
as a run on $\Board_n$. Replies inside $\Board_{17}$ are read
verbatim. A reply beyond $\Board_{17}$ is
replaced, in the run on $\Board_{23}$, by a \emph{stand-in}: a fixed cell beyond
$\Board_{17}$ when the pair has no hidden corner, and an
\emph{anchor}, a cell beyond $\Board_{17}$ that forms with the other
reply the same hidden corner, when it has one. The \emph{anchor board}
$\Board_{23}$ is the least board holding every anchor
(Remark~\ref{rem:anchors}).

A \emph{history} is the sequence of moves played so far, and for a
history $h$ and a sequence $t$ of further moves $h\,t$ is their
concatenation. A strategy $\sigma$ is given the board, and
$\sigma_n(h)$ is its move after $h$ on $\Board_n$. The guards are those
of Section~\ref{sec:verify}.

\begin{theorem}\label{thm:simulate}
Let $\sigma$ be a strategy, let $m, n \ge 17$, and let
$x_1, x_2 \in \Board_6$. Let $h = (x_1, o_1, x_2, o_2)$ be a history
on $\Board_n$ and $h' = (x_1, o_1', x_2, o_2')$ a history on
$\Board_m$ such that for each $i$ either $o_i' = o_i \in \Board_{17}$,
or $o_i \in \Board_n \setminus \Board_{17}$ and
$o_i' \in \Board_m \setminus \Board_{12}$. Let $g$ be either no cell or
a cell of $\Board_6$. Suppose that
\begin{enumerate}
\item[(i)] on every history extending $h$ in which the first player's
  stones lie in $\Board_6$, $\sigma$ plays a cell completing one of the
  first player's squares whenever one exists;
\item[(ii)] $\sigma_n(h\,t) = \sigma_m(h'\,t)$ for every sequence $t$
  of further moves;
\item[(iii)] the forced tree of $\sigma$ on $\Board_m$ below $h'$
  passes the guards, with $g$ as the cell of (G4) if $g$ is a cell, and
  every branch ends in a first-player square;
\item[(iv)] if $o_1$ or $o_2$ lies beyond $\Board_{17}$, then either
  $g$ is no cell and $\{o_1, o_2\}$ has no hidden corner
  (Definition~\ref{def:latent}), in particular if one of them has no
  coordinate at least $12$, or $g$ is the hidden corner of
  $\{o_1, o_2\}$.
\end{enumerate}
Then the tree in which the first player plays $\sigma$'s move at every
node below $h$ on $\Board_n$ is winning.
\end{theorem}

\begin{proof}
By induction along the forced tree of (iii). Consider a node of that
tree together with the corresponding real position on $\Board_n$; the
two differ only at the two reply slots, and the blocks played so far,
all in $\Board_6$ by (G2), are shared. The two positions agree inside
$\Board_{12}$: a real reply beyond $\Board_{17}$ and its stand-in both
lie outside $\Board_{12}$, and every other stone is shared. By (ii)
$\sigma$ plays the same cell $m$ in both. If $m$ completes a square in
the anchor position, the square has three corners among the first
player's stones, which lie in $\Board_6$, so it lies in $\Board_{12}$
by Lemma~\ref{lem:span}, and $m$ is empty and completes it on
$\Board_n$ too. Otherwise $m \in \Board_6$ by (G2), and the threat
cells of $m$ are the same on both boards, since each threat square has
three corners in $\Board_6$ and so lies in $\Board_{12}$; by (G3) the
second player has no threat in the anchor position.

No reply $r$ completes a second-player square on $\Board_n$. Suppose
$\Sq$ were such a square.

Case 1: $\Sq \subseteq \Board_{17}$. Its corners other than $r$ are
second-player stones inside $\Board_{17}$, so blocks or replies read
verbatim, all shared with the anchor position, since the real replies
beyond $\Board_{17}$ lie outside $\Sq$ and a stand-in only adds a
stone. Then either $r$ is a stand-in and the anchor position already
holds $\Sq$, or $r$ is empty in the anchor position and completes
$\Sq$ there. The second case contradicts (G3) at the current node. In
the first case, the anchor position holds no second-player square: the
two replies alone are too few corners, so the last corner placed was a
block, and that block answered a threat that was standing after the
previous move of $\sigma$, again contradicting (G3).

Case 2: $\Sq \not\subseteq \Board_{17}$. Two blocks among its corners
would confine it to $\Board_{12}$ by Lemma~\ref{lem:span}, so
its corners are $r$, one block $d \in \Board_6$, and both replies, and
a reply lies beyond $\Board_{17}$: any three corners of a square
include a diagonally opposite pair, and such a pair in $\Board_{17}$
confines $\Sq$ there by Lemma~\ref{lem:span}. By
Lemma~\ref{lem:dichotomy} $d$ is the
hidden corner of $\{o_1, o_2\}$, which contradicts (iv): either the
pair has none, or its hidden corner is $g$, and $g$ is never a block,
because a block at $g$ arises only when $g$ is the only threat cell,
which (G4) excludes, or when the second player chooses $g$ among two
threat cells, in which case $\sigma$ completes the other by (i) and
the branch ends.

A reply that blocks nothing leaves a threat, which $\sigma$ completes
by (i). A block $r \ne g$ is a block in the anchor run as well, and the
induction continues along the anchor branch, whose invariants are
preserved: the new block lies in $\Board_6$ and is not $g$.
\end{proof}

The formal theorem is stated for any node of the anchor run reached by
forced play, with the blocks so far in $\Board_6$ and $g$ not among
the second player's stones, and is proved by induction on the checker's
fuel.

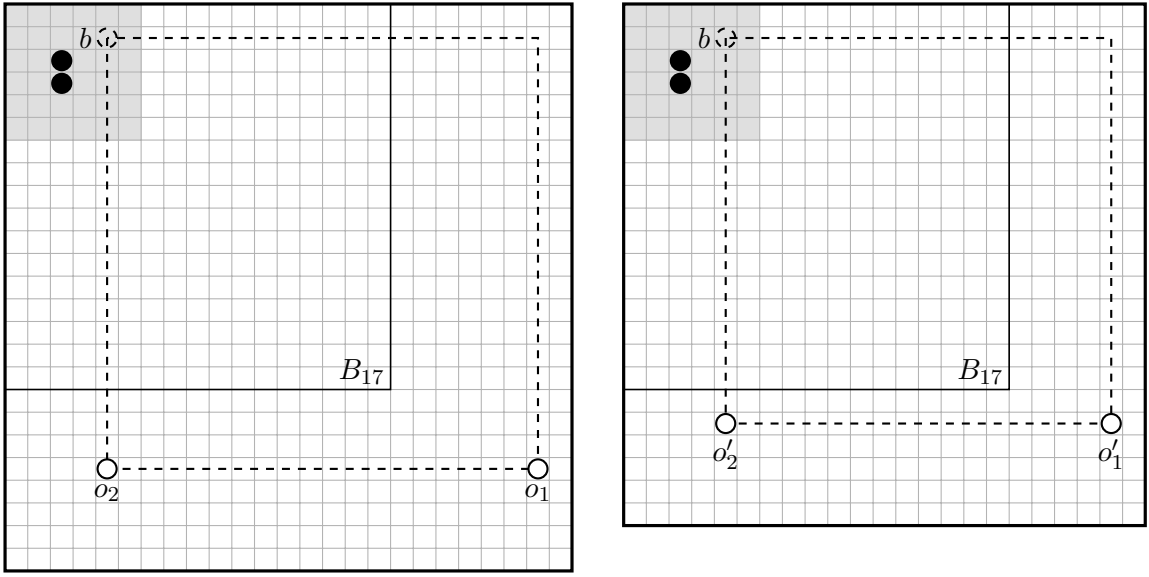
\begin{figure}[ht]
\centering
% Generated by python3 -m squares.paper.figures; do not edit.
\begin{tikzpicture}[scale=0.3,baseline=(current bounding box.north)]
  \fill[gray!25] (0,0) rectangle (6,-6);
  \draw[gray!60] (0,0) grid (25,-25);
  \draw[semithick] (0,0) rectangle (17,-17);
  \draw[very thick] (0,0) rectangle (25,-25);
  \draw[thick,dashed] (4.5,-1.5) rectangle (23.5,-20.5);
  \foreach \r/\c in {2/2, 3/2}
    \filldraw[fill=black,thick] (\c+0.5,-\r-0.5) circle (0.42);
  \foreach \r/\c in {20/23, 20/4}
    \filldraw[fill=white,thick] (\c+0.5,-\r-0.5) circle (0.42);
  \draw[thick,dashed] (4.5,-1.5) circle (0.42);
  \node[anchor=east,inner sep=1pt] at (3.96,-1.5) {$b$};
  \node[anchor=north,inner sep=1pt] at (23.5,-21.04) {$o_1$};
  \node[anchor=north,inner sep=1pt] at (4.5,-21.04) {$o_2$};
  \node[anchor=south east,inner sep=2pt] at (17,-17) {$\Board_{17}$};
\end{tikzpicture}
\hspace{1em}
\begin{tikzpicture}[scale=0.3,baseline=(current bounding box.north)]
  \fill[gray!25] (0,0) rectangle (6,-6);
  \draw[gray!60] (0,0) grid (23,-23);
  \draw[semithick] (0,0) rectangle (17,-17);
  \draw[very thick] (0,0) rectangle (23,-23);
  \draw[thick,dashed] (4.5,-1.5) rectangle (21.5,-18.5);
  \foreach \r/\c in {2/2, 3/2}
    \filldraw[fill=black,thick] (\c+0.5,-\r-0.5) circle (0.42);
  \foreach \r/\c in {18/21, 18/4}
    \filldraw[fill=white,thick] (\c+0.5,-\r-0.5) circle (0.42);
  \draw[thick,dashed] (4.5,-1.5) circle (0.42);
  \node[anchor=east,inner sep=1pt] at (3.96,-1.5) {$b$};
  \node[anchor=north,inner sep=1pt] at (21.5,-19.04) {$o_1'$};
  \node[anchor=north,inner sep=1pt] at (4.5,-19.04) {$o_2'$};
  \node[anchor=south east,inner sep=2pt] at (17,-17) {$\Board_{17}$};
\end{tikzpicture}
\caption{Left: a real history on $\Board_{25}$, with $\Board_{17}$
outlined and $\Board_6$ shaded. The first player holds $(2,2)$
and $(3,2)$; the replies $o_1 = (20,23)$ and $o_2 = (20,4)$, both
beyond $\Board_{17}$, are two corners of the dashed square
$\Sq(1,4,19)$ with hidden corner $b = (1,4)$. Right: the history on
the anchor board $\Board_{23}$; the anchors $o_1' = (18,21)$ and
$o_2' = (18,4)$ are two corners of $\Sq(1,4,17)$ with the same hidden
corner. The two histories have the same abstraction, so the rule plays
the same line from each.}
\label{fig:agreement}
\end{figure}

For the rule, hypotheses (i) and (ii) are lemmas about its definition.
Hypothesis (i) is the hypothesis of Proposition~\ref{prop:forced},
which the rule satisfies as noted there. For
(ii), the rule's move depends on a history only through its
abstraction and through the second move of the first reply. A reply
beyond $\Board_{17}$ and its stand-in both drop out of the abstraction,
a reply inside $\Board_{17}$ is its own stand-in, and the stand-in pair
is chosen with the same hidden corner as the real pair, so the two
abstractions agree at every node; the second moves agree because every
first reply beyond $\Board_6$, and every stand-in, has the default
second move $(3,2)$. Figure~\ref{fig:agreement} shows a real history
and its anchor history.

Pairs with a hidden corner $b = (b_1, b_2)$ are covered by three
families of \emph{anchored runs}, each a forced run of the rule on
$\Board_{23}$ from the history $((2,2), o_1', x_2(o_1'), o_2')$ with
$b$ as the cell of (G4), preceded by a check that $b$ is the hidden
corner of the anchor pair.
\begin{itemize}
\item The \emph{far} family, for both replies beyond $\Board_{17}$:
  $o_1' = (b_1 + 17, b_2 + 17)$ and $o_2' = (b_1 + 17, b_2)$, two
  corners of $\Sq(b_1, b_2, 17)$; one run per $b \in \Board_6$, $36$
  runs.
\item The \emph{exact} family, for $o_1 \in \Board_{17}$ and $o_2$
  beyond it. By Lemma~\ref{lem:dichotomy}, the square through $b$ and
  the pair is $\Sq(b_1, b_2, s)$ with $s \ge 12$, and it leaves
  $\Board_{17}$. The corner opposite $b$ bounds the square, so $o_1$ is
  not that corner; hence $o_1$ lies in the row or the column of $b$ at
  distance $s$, and $o_1 \in \Board_{17}$ puts that coordinate of $o_1$
  below $17$. Say $o_1 = (a, w)$ and $b = (a, h)$ with $a, h < 6$,
  $12 \le w < 17$, and $s = w - h$. Then $o_1' = o_1$ and
  $o_2' = (a + s, w)$, the corner of $\Sq(a, h, s)$ opposite $b$, which
  lies beyond $\Board_{17}$ exactly when the square does. One run per
  $(a, w, h)$ with $a + w - h \ge 17$ and one per transposed triple,
  $70$ runs.
\item The \emph{far-near} family, for $o_1$ beyond $\Board_{17}$ and
  $o_2 \in \Board_{17}$: the exact family with the two replies
  exchanged, $70$ runs.
\end{itemize}
This gives $\anchoredPositions$ anchored runs in all.

\begin{proof}[Proof of Theorem~\ref{thm:uniform}]
Let $n \ge 17$. The rule opens at $(2,2)$ and answers the first reply
$o_1$ with $x_2 = x_2(o_1) \in \Board_6$; both moves are legal and
complete nothing, so it remains to show that the rule's tree below
$h = ((2,2), o_1, x_2, o_2)$ is winning for every reply pair. Apply
Theorem~\ref{thm:simulate} with $m = 23$, $\sigma$ the rule, and $h'$,
$g$ as follows, calling a reply \emph{far} if it lies beyond
$\Board_{17}$.
\begin{enumerate}
\item[(a)] Neither reply far: $h' = h$, $g$ no cell, and (iii) is the
  branch at $h$ of the forced run of the rule on $\Board_{23}$; (iv) is
  vacuous.
\item[(b)] A far reply and no hidden corner: a far $o_1$ is replaced by
  $(22,22)$ and a far $o_2$ by $(22,21)$, $g$ is no cell, and (iii) is
  again a branch of the run on $\Board_{23}$. A stand-in forms no hidden
  corner with a cell of $\Board_{17}$: it shares no row or column with
  such a cell, and if the two are diagonally opposite corners of a
  square, the other two corners have a row or a column of at least
  $21$. The two stand-ins share row $22$, so a hidden corner of theirs
  would lie in row $21$ or $23$. So the stand-in pair has no hidden
  corner either.
\item[(c)] Both replies far with hidden corner $b$: the far family,
  $g = b$.
\item[(d)] $o_1 \in \Board_{17}$, $o_2$ far, hidden corner $b$: the
  exact family, $g = b$.
\item[(e)] $o_1$ far, $o_2 \in \Board_{17}$, hidden corner $b$: the
  far-near family, $g = b$.
\end{enumerate}
In (c) to (e), (iii) is the anchored run, and (iv) holds by the case
hypothesis. Hypothesis (i) is the lemma above; (ii) holds because the
real and anchor pairs have the same hidden corner, checked at the head of
each anchored run and shown in (b) for the stand-ins, and the same
second move, so that the abstractions of $h\,t$ and $h'\,t$ agree for
every $t$.
\end{proof}

\begin{remark}\label{rem:anchors}
The anchors of the far family reach coordinate $5 + 17 = 22$, those of
the other two families $a + s \le 5 + 16 = 21$, and the stand-ins sit
at $22$, so $\Board_{23}$ is the least board holding them all; any
larger board would serve. The anchor board, the stand-ins, and the
anchors are choices made in the proof, and the rule does not depend on
them.
\end{remark}

\begin{proof}[Proof of Theorems~\ref{thm:main} and~\ref{thm:rule}]
Theorem~\ref{thm:rule} follows from the checked runs at $6$ to $17$
with Proposition~\ref{prop:forced}, and from Theorem~\ref{thm:uniform}
for $n \ge 17$. Theorem~\ref{thm:main} follows from
Propositions~\ref{prop:small} and~\ref{prop:five} and from
Theorem~\ref{thm:rule} with the moves forgotten.
\end{proof}

\section{Formalization}\label{sec:lean}

The Lean development defines cells, boards, and squares as in
Section~\ref{sec:defs} (\texttt{Core/Square.lean}), and the winning
predicate \texttt{P1CanWin} (\texttt{Game/Rules.lean}) as an inductive
family whose derivations are exactly strategy trees: one first-player
move per node, every reply covered, and no derivation for a game that
ends in a draw. The second player's win is defined from it by swapping
roles after a move; the first player's ability to secure a draw is a
second inductive family, whose branches may also end on a full board,
and the second player's is defined from that by the same swap
(\texttt{Game/Outcome.lean}). Strategy stealing
(\texttt{second\_player\_no\_win}, with the monotonicity step
\texttt{P1CanWin.insert\_left}) and determinacy
(\texttt{P1CanWin.or\_p2CanDraw}) are proved against these definitions
with no computation. These definitions, a few dozen lines, are the
whole specification; everything else is proved against them, and the
development has no \texttt{sorry} and no unproved declaration. The two
main statements are
\begin{center}\small
\texttt{outcome (n : $\mathbb{N}$) : (P1CanWin n $\varnothing$
$\varnothing$ $\leftrightarrow$ 5 $\le$ n) $\wedge$ ((P1CanDraw n
$\varnothing$ $\varnothing$ $\wedge$ P2CanDraw n $\varnothing$
$\varnothing$) $\leftrightarrow$ n $\le$ 4)}\\[2pt]
\texttt{plus\_wins\_with \char123 n : $\mathbb{N}$\char125\ (h : 6 $\le$ n) :
WinsWith plusMove n []}
\end{center}
where \texttt{P1CanWin n A B} says that the player to move, holding
$A$ against $B$ on $\Board_n$, can force a win. Running
\texttt{lake build} checks everything, \texttt{lake exe checker n}
replays the certificate for $\Board_n$, and \texttt{\#print axioms}
lists the axioms of a theorem.

The predicate behind Theorem~\ref{thm:rule} is \texttt{WinsWith}. For a
strategy $\sigma$, a board size $n$, and a history $h$ with the first
player to move, $\texttt{WinsWith}\ \sigma\ n\ h$ holds if $\sigma$'s
move at $h$ is legal and either completes a square or admits at least
one reply, no reply completes a second-player square, and
\texttt{WinsWith} holds again after every reply. It has the
constructors of
\texttt{P1CanWin} with one extra premise on each, that the move is the
one $\sigma$ returns, and implies \texttt{P1CanWin} by forgetting that
premise. A checker that plays $\sigma$ at every node proves
\texttt{WinsWith} only if its soundness proof never invents a move,
which is why Proposition~\ref{prop:forced} needs its hypothesis,
\texttt{TakesWinsOn}.

Three executable checkers are proved sound: \texttt{check}, which
replays a table against all reply sequences (\texttt{check\_sound});
the forced checker, generic over strategies
(\texttt{checkForcedS\_winsWith}, Proposition~\ref{prop:forced}); and
\texttt{refute} (\texttt{refute\_sound}, Proposition~\ref{prop:small}).
Each proof is an induction over the checker's recursion, whose guards
supply the premises of the constructors, so soundness never inspects
the strategy. The $5 \times 5$ table is
\texttt{fiveTable} and Proposition~\ref{prop:five} is
\texttt{p1\_wins\_five}.

The rule is the definition \texttt{plusMove}; the file
\texttt{Plus/Rule.lean}, which defines it with the abstraction and the
line search, is about $\ruleLeanLines$ lines.
Hypothesis (i) of Theorem~\ref{thm:simulate} for the rule is
\texttt{plusMove\_takesWins} and hypothesis (ii) is
\texttt{plusMove\_agree}; both are proved by unfolding the definition,
and neither is an evaluation. The simulation theorem is \texttt{simulate}, and
Theorem~\ref{thm:uniform} is
\texttt{plus\_wins\_with\_of\_seventeen\_le}, whose body
\texttt{plus\_uniform\_aux} dispatches the cases (a) to (e) to
\texttt{simulate} with the run on the anchor board
(\texttt{plus\_forced\_anchor}) and the three families
\texttt{farAnchorsSafe}, \texttt{exactAnchorsSafe}, and
\texttt{farNearAnchorsSafe}. Lemma~\ref{lem:span} is not a separate
Lean theorem; the containments it gives are proved where they are
used: \texttt{three\_corners\_confine} and
\texttt{square\_filter\_confine} confine a threat square to the
region, and \texttt{escape\_pick} and \texttt{escape\_farm} confine a
square with two region corners. The hidden corner of a pair
is computed by \texttt{farmKey}, which the abstraction and the anchored
runs use; the Lean identifiers call the hidden corner a \emph{farm}
(\texttt{farmKey}, \texttt{Key.farm}) and the cell of (G4) \texttt{hot}.
A cell \texttt{farmKey} returns lies in $\Board_6$ and is a corner of a
square through the pair (\texttt{farmKey\_corners}), and the second
claim of Lemma~\ref{lem:dichotomy} is \texttt{escape\_farm}: for a
square as in the lemma, \texttt{farmKey} returns $b$. Hypothesis (iv)
of Theorem~\ref{thm:simulate} is stated through \texttt{farmKey}: $g$
is no cell and \texttt{farmKey} returns none, or \texttt{farmKey}
returns $g$. The uniqueness
claim is not needed by the formal proof, which reads the one value
\texttt{farmKey} returns. Lemma~\ref{lem:arms} is not used by the formal proof,
which checks the lines directly. Theorem~\ref{thm:rule} is
\texttt{plus\_wins\_with} and Theorem~\ref{thm:main} is
\texttt{outcome}, whose first half is \texttt{p1\_wins\_iff}.

\begin{table}[ht]
\centering
\begin{tabular}{rrrcrrr}
\toprule
$n$ & pairs & time (s) & & $n$ & pairs & time (s) \\
\midrule
5 & full tree & \timeFive & & 12 & \pairsTwelve & \timeTwelve \\
6 & \pairsSix & \timeSix & & 13 & \pairsThirteen & \timeThirteen \\
7 & \pairsSeven & \timeSeven & & 14 & \pairsFourteen & \timeFourteen \\
8 & \pairsEight & \timeEight & & 15 & \pairsFifteen & \timeFifteen \\
9 & \pairsNine & \timeNine & & 16 & \pairsSixteen & \timeSixteen \\
10 & \pairsTen & \timeTen & & 17 & \pairsSeventeen & \timeSeventeen \\
11 & \pairsEleven & \timeEleven & & 23 & \pairsAnchor & \timeAnchor \\
\bottomrule
\end{tabular}
\medskip
\caption{The runs of the checkers, timed with the checker executable
on a MacBook Pro (16-inch, 2024) with an Apple M4 Max chip and 36 GB of
memory. At $5$ the enumerating checker
walks the full tree of the table; from $6$ up the forced checker runs
the rule, with the number of reply pairs at the root. The run at $23$
is the run on the anchor board; the $\anchoredPositions$ anchored runs
take $\timeAnchored$ seconds together.}
\label{tab:runs}
\end{table}

Concrete boards are discharged by native evaluation
(\texttt{native\_decide}): the refutations at $1$ to $4$, the
$5 \times 5$ run, one evaluation for the rule's runs at sizes $6$ to
$16$, one for $17$, one for the anchor board, and one per anchored
family. The axioms of \texttt{outcome} and \texttt{plus\_wins\_with} are
therefore \texttt{propext}, \texttt{Classical.choice},
\texttt{Quot.sound}, and one instance of Lean's evaluation axiom per
certificate. The kernel checks the soundness theorems; the evaluations
are performed by the Lean compiler and runtime, which are therefore
trusted, the usual trust base of computer-assisted proofs in Lean. The
whole development builds in about $\buildMinutes$ minutes on the
machine of Table~\ref{tab:runs}, almost all of it the evaluations.

The Lean development and the Python code, with the $5 \times 5$
solver, the implementation of the rule, and the scripts that produce
the numbers and board diagrams of this paper, accompany the paper as
ancillary material and are at
\url{https://github.com/sean-reid/square-achievement}, tag
\texttt{v1.0.0}. The development builds with Lean~4.33.1 and mathlib at
its tag \texttt{v4.33.1}.

\section{Open questions}

The rule wins by ply $\latestPly$. What is the least ply by which the
first player can force a win on $\Board_n$, and is there a rule of
comparable length that attains it on every board? Other orderings of
the lines, or other target shapes, might shorten the longest games at
the cost of a longer statement.

The second move has seven exceptions. Is there an opening, or an
opening with a fixed second move, from which a rule of this kind wins
with none? The search of Section~\ref{sec:rule} did not vary the
opening.

Whether a short rule exists for $\Board_5$, where the table needs quiet
moves, is open.

Base $16$ would serve Lemma~\ref{lem:dichotomy}
(Remark~\ref{rem:base}); whether the rule still wins when it reads only
$\Board_{16}$ was not checked. Whether a rule of this kind can play
inside a region smaller than $\Board_6$, or read less than the second
player's stones inside its base board, for instance only those inside
$\Board_6$ with a summary of the rest, is open.

\section*{Acknowledgements}

I thank Nicholas Tan for many games and for productive discussions that
shaped the strategy presented here, and Anthony Gabriele for reviewing
the presentation and content of the paper.

\bibliographystyle{amsplain}
\bibliography{refs}

\end{document}